\documentclass[a4paper,12pt]{article}

\usepackage[latin1]{inputenc}
\usepackage[francais,english]{babel}

\usepackage[T1]{fontenc}
\usepackage[scaled]{helvet}

\usepackage{fullpage}
\usepackage{url}
\usepackage{a4,amsmath,amssymb,amsfonts,amsthm}
\usepackage{amsrefs}

\usepackage[small,nohug,heads=vee]{diagrams} 
\diagramstyle[labelstyle=\scriptstyle]

\usepackage[all,cmtip]{xy}

\begin{document}


\newcommand{\fr}{\selectlanguage{francais}}
\newcommand{\en}{\selectlanguage{english}}


\def\M#1{\mathbb#1} 
\def\mR{\M{R}}           
\def\mZ{\M{Z}}           
\def\mN{\M{N}}           
\def\mQ{\M{Q}}       
\def\mC{\M{C}}  
\def\mG{\M{G}}
\def\mP{\M{P}}


\def\Spec{{\rm Spec}}
\def\rg{{\rm rg}}
\def\Hom{{\rm Hom}}
\def\Aut{{\rm Aut}}
 \def\Tr{{\rm Tr}}
 \def\Exp{{\rm Exp}}
 \def\Gal{{\rm Gal}}
 \def\End{{\rm End}}
 \def\det{{{\rm det}}}
 \def\Td{{\rm Td}}
 \def\ch{{\rm ch}}
 \def\che{{\rm ch}_{\rm eq}}
  \def\Spec{{\rm Spec}}
\def\Id{{\rm Id}}
\def\Zar{{\rm Zar}}
\def\Supp{{\rm Supp}}
\def\eq{{\rm eq}}
\def\Ann{{\rm Ann}}
\def\LT{{\rm LT}}
\def\Pic{{\rm Pic}}
\def\rg{{\rm rg}}
\def\et{{\rm et}}
\def\sep{{\rm sep}}
\def\ppcm{{\rm ppcm}}
\def\ord{{\rm ord}}
\def\Gr{{\rm Gr}}
\def\rk{{\rm rk}}
\def\Stab{{\rm Stab}}
\def\im{{\rm im}}
\def\Sm{{\rm Sm}}
\def\red{{\rm red}}
\def\Frob{{\rm Frob}}
\def\Ver{{\rm Ver}}
\def\Div{{\rm Div}}
\def\ab{{\rm ab}}


\def\beginProof{\par{\bf Proof. }}
 \def\endProof{${\qed}$\par\smallskip}
 \def\pr{^{\prime}}
 \def\prpr{^{\prime\prime}}
 \def\mtr#1{\overline{#1}}
 \def\ra{\rightarrow}
 \def\mfp{{\mathfrak p}}
 \def\mfm{{\mathfrak m}}
 
 \def\char{{\rm char}}
  \def\mQ{{\Bbb Q}}
 \def\mR{{\Bbb R}}
 \def\mZ{{\Bbb Z}}
 \def\mC{{\Bbb C}}
 \def\mN{{\Bbb N}}
 \def\mF{{\Bbb F}}
 \def\mA{{\Bbb A}}
  \def\mG{{\Bbb G}}
 \def\CI{{\cal I}}
 \def\CA{{\cal A}}
 \def\CD{{\cal D}}
 \def\CE{{\cal E}}
 \def\CJ{{\cal J}}
 \def\CH{{\cal H}}
 \def\CO{{\cal O}}
 \def\CA{{\cal A}}
 \def\CB{{\cal B}}
 \def\CC{{\cal C}}
 \def\CK{{\cal K}}
 \def\CL{{\cal L}}
 \def\CI{{\cal I}}
 \def\CM{{\cal M}}
  \def\CN{{\cal N}}
\def\CP{{\cal P}}
 \def\CZ{{\cal Z}}
\def\CR{{\cal R}}
\def\CG{{\cal G}}
\def\CX{{\cal X}}
\def\CY{{\cal Y}}
\def\CV{{\cal V}}
\def\CW{{\cal W}}
\def\CT{{\cal T}}
 \def\wt#1{\widetilde{#1}}
 \def\mod{{\rm mod\ }}
 \def\refeq#1{(\ref{#1})}
 \def\blb{{\big(}}
 \def\brb{{\big)}}
\def\mc{{{\mathfrak c}}}
\def\mcpr{{{\mathfrak c}'}}
\def\mcprpr{{{\mathfrak c}''}}
\def\ss{{\rm ss}}
\def\parf{{\rm parf}}
\def\P1{{{\bf P}^1}}
\def\cod{{\rm cod}}
\def\pr{\prime}
\def\prpr{\prime\prime}
\def\ss{\scriptstyle}
\def\OX{{ {\cal O}_X}}
\def\mpartial{{\mtr{\partial}}}
\def\inv{{\rm inv}}
\def\indlim{\underrightarrow{\lim}}
\def\prolim{\underleftarrow{\lim}}
\def\pprolim{'\prolim'}
\def\Pro{{\rm Pro}}
\def\Ind{{\rm Ind}}
\def\Ens{{\rm Ens}}
\def\without{\backslash}
\def\pbdb{{\Pro_b\ D^-_c}}
\def\qc{{\rm qc}}
\def\Com{{\rm Com}}
\def\an{{\rm an}}
\def\gfield{{\rm\bf k}}
\def\s{{\rm s}}
\def\dR{{\rm dR}}
\def\ari#1{\widehat{#1}}
\def\ul#1{\underline{#1}}
\def\sul#1{\underline{\scriptsize #1}}
\def\mou{{\mathfrak u}}
\def\ich{\mathfrak{ch}}
\def\cl{{\rm cl}}
\def\K{{\rm K}}
\def\R{{\rm R}}
\def\F{{\rm F}}
\def\L{{\rm L}}
\def\pgcd{{\rm pgcd}}
\def\rc{{\rm c}}
\def\N{{\rm N}}
\def\E{{\rm E}}
\def\H{{\rm H}}
\def\CHOW{{\rm CH}}
\def\A{{\rm A}}
\def\d{{\rm d}}
\def\Res{{\rm  Res}}
\def\GL{{\rm GL}}
\def\Alb{{\rm Alb}}
\def\alb{{\rm alb}}
\def\Hdg{{\rm Hdg}}
\def\Num{{\rm Num}}
\def\Irr{{\rm Irr}}
\def\Frac{{\rm Frac}}
\def\Sym{{\rm Sym}}
\def\TV{{\cal Z}}
\def\indlim{\underrightarrow{\lim}}
\def\prolim{\underleftarrow{\lim}}
\def\LE{{\rm LE}}
\def\LEG{{\rm LEG}}
\def\MFT{{\mathfrak T}}
\def\Exc{{\rm Exc}}
\def\Crit{{\rm Crit}}
\def\MFW{{\mathfrak W}}
\def\Transp{{\rm Transp}}


\def\RHom{{\rm RHom}}
\def\rRHom{{\mathcal RHom}}
\def\rHom{{\mathcal Hom}}
\def\dotimes{{\overline{\otimes}}} 
\def\Ext{{\rm Ext}}
\def\rExt{{\mathcal Ext}}
\def\Tor{{\rm Tor}}
\def\rTor{{\mathcal Tor}}
\def\SP{{\mathfrak S}}
\def\WR{{\mathfrak R}}

\def\H{{\rm H}}
\def\D{{\rm D}}
\def\Del{{\mathfrak D}}
\def\alg{{\rm alg}}
\def\sh{{\rm sh}}
 

 \newtheorem{theor}{Theorem}[section]
 \newtheorem{prop}[theor]{Proposition}
 \newtheorem{propdef}[theor]{Proposition-Definition}
 \newtheorem{sublemma}[theor]{sublemma}
 \newtheorem{cor}[theor]{Corollary}
 \newtheorem{lemma}[theor]{Lemma}
  \newtheorem{lemme}[theor]{Lemme}
 \newtheorem{sublem}[theor]{sub-lemma}
 \newtheorem{defin}[theor]{Definition}
 \newtheorem{conj}[theor]{Conjecture}
 \newtheorem{quest}[theor]{Question}
  \newtheorem{remark}[theor]{Remark}


\author{
 Damian R\"OSSLER\footnote{Institut de Math\'ematiques, 
Equipe Emile Picard, 
Universit\'e Paul Sabatier, 
118 Route de Narbonne,  
31062 Toulouse cedex 9, 
FRANCE, E-mail: rossler@math.univ-toulouse.fr}}

\title{On the Manin-Mumford and Mordell-Lang conjectures in positive characteristic}
\maketitle

\abstract{We prove that in positive characteristic, the Manin-Mumford conjecture implies the Mordell-Lang conjecture, in the situation where the ambient variety 
is an abelian variety defined over the function field of a smooth curve over a finite field and 
the relevant group is a finitely generated group. In particular, in the setting of 
the last sentence, we provide a proof of the Mordell-Lang conjecture, which does not depend on tools coming from model theory.}

\parskip=3pt
\parindent=0pt

\section{Introduction}

Let $B$ be a semiabelian variety over an algebraically closed field $F$ of 
characteristic $p>0$. Let $Y$ be an irreducible reduced closed subscheme of $B$. 
Let $\Lambda\subseteq B(F)$ be a subgroup. Suppose that $\Lambda\otimes_\mZ\mZ_{(p)}$ is 
a finitely generated $\mZ_{(p)}$-module (here, as is customary, we write $\mZ_{(p)}$ for the localization 
of $\mZ$ at the prime $p$). 

Let $C:=\Stab(Y)^\red$, where $\Stab(Y)=\Stab_B(Y)$ is the translation stabilizer of $Y$. 
This is the closed subgroup scheme  of $B$, which is characterized uniquely 
by the fact that for any scheme $S$ and any morphism $b: S\to B$, translation by $b$ on the
product $B\times_F S$ maps the subscheme $Y\times S$ to itself if and
only if $b$ factors through $\Stab_B(Y)$.  Its existence is proven in 
\cite[exp.~VIII, Ex.~6.5~(e)]{SGA3-2}.

The Mordell-Lang conjecture for $Y$ and $B$ is now the following statement. 

\begin{theor}[Mordell-Lang conjecture; Hrushovski \cite{Hrushovski-Mordell-Lang}]
If $Y\cap\Lambda$ is Zariski dense in $Y$ then there is 
\begin{itemize}
\item a semiabelian variety $B'$ over $F$; 
\item a homomorphism with finite kernel $h: B'\to B/C$;
\item a model $\bf B'$ of $B'$ over a finite subfield $\mF_{p^r} \subset F$;
\item an irreducible reduced closed subscheme ${\bf Y'}\hookrightarrow{\bf B'}$;
\item a point $b\in (B/C)(F)$, such that\,\, 
$Y/C = b + h_*\bigl({\bf Y'}\times_{\mF_{p^r}}F\bigr).$
\end{itemize}
\label{MLtheor}
\end{theor}

Here $h_*\bigl({\bf Y'}\times_{\mF_{p^r}}F\bigr)$ refers to the scheme-theoretic image of 
${\bf Y'}\times_{\mF_{p^r}}F$ by $h$. Since $h$ is finite and ${\bf Y'}\times_{\mF_{p^r}}F$ is reduced, this implies that $h_*\bigl({\bf Y'}\times_{\mF_{p^r}}F\bigr)$ is simply the set-theoretic image of ${\bf Y'}\times_{\mF_{p^r}}F$ by $h$, endowed with its reduced-induced scheme structure.

Theorem \ref{MLtheor} in particular implies the following result, 
which will perhaps seem more striking on first reading. Suppose that there are no non-trivial homomorphisms from $B$ to a 
semiabelian variety, which has a model over a finite field. Then : 
if $Y\cap\Lambda$ is Zariski dense in $Y$ then $Y$ is the translate of 
an abelian subvariety of $B$. 

Theorem \ref{MLtheor} was first proven in 1996  by E. Hrushovski using deep results from model theory, in particular the Hrushovski-Zilber theory of Zariski geometries (see \cite{Hrushovski-Zilber-Zariski}).  An algebraic proof of Theorem \ref{MLtheor} in the situation 
where $B$ is an ordinary abelian variety was given by D. Abramovich and J.-.F Voloch in \cite{Voloch-Towards}.
In the situation where $Y$ is a smooth curve embedded into $B$ as 
its Jacobian, the theorem was known to be true much earlier. See for instance \cite{Samuel-Complements} 
and \cite{Szpiro-Seminaire-Pinceaux}. The earlier proofs for curves relied on the use of heights, which 
do not appear in the later approach of Voloch and Hrushovski, which 
is parallel and inspired by A. Buium's approach in characteristic $0$ via 
differential equations (see below). 

The {\it Manin-Mumford conjecture} has exactly the same form as the Mordell-Lang conjecture, 
but $\Lambda$ is replaced by the group 
$\Tor(B(F))$ of points of finite order of $B(F)$. For the record, we state it in full: 

\begin{theor}[Manin-Mumford conjecture; Pink-R\"ossler \cite{PR2}]
If $Y\cap\Tor(B(F))$ is Zariski dense in $Y$ then there is 
\begin{itemize}
\item a semiabelian variety $B'$ over $F$; 
\item a homomorphism with finite kernel $h: B'\to B/C$;
\item a model $\bf B'$ of $B'$ over a finite subfield $\mF_{p^r} \subset F$;
\item an irreducible reduced closed subscheme ${\bf Y'}\hookrightarrow{\bf B'}$;
\item a point $b\in (B/C)(F)$, such that\,\, 
$Y/C = b + h_*\bigl({\bf Y'}\times_{\mF_{p^r}}F\bigr).$
\end{itemize}
\label{MMtheor}
\end{theor}

See also \cite{Scanlon-A-positive} for a model-theoretic proof of the Manin-Mumford conjecture. 

{\bf Remark (important).} Notice that the Manin-Mumford conjecture is {\it not} a special case of the Mordell-Lang conjecture, because $\Tor(A(F))$ is not in general a finitely generated $\mZ_{(p)}$-module (because
$\Tor(A(F))[p^\infty]$ is not finite in general). Nevertheless, it seems reasonable to conjecture that 
Theorem \ref{MLtheor}  should still be true when the the hypothesis that $\Lambda\otimes_\mZ\mZ_{(p)}$ is finitely 
generated is 
replaced by the weaker hypothesis that $\Lambda\otimes_\mZ\mQ$ is finitely generated. 
This last statement, which is still not proven in general, is often called the {\it full Mordell-Lang conjecture} and it would have Theorems \ref{MLtheor} and \ref{MMtheor} as special cases. See 
\cite{Ghioca-Moosa-Division} for more about this.

Now suppose that the group $\Lambda$ is actually finitely generated and that $B$ arises by base-change to 
$F$ from an abelian variety $B_0$, which is defined over a function field of transcendence degree $1$ over a finite field. The main result of this text is then the proof of the fact that the 
Manin-Mumford conjecture in general implies the Mordell-Lang conjecture in this situation. 
We follow here the lead of A. Pillay, who suggested in a talk he gave in Paris on Dec. 
17th 2010 that it should be possible to establish this logical link without proving 
the Mordell-Lang conjecture first. See Theorem \ref{MMprop} and its corollary below for a precise statement. 

The interest of an algebro-geometric (in contrast with model-theoretic) proof of the implication Manin-Mumford $\Longrightarrow$ Mordell-Lang is that it provides in particular 
an algebro-geometric proof of the Mordell-Lang conjecture.

Let $K_0$ be the function field of a smooth curve over $\bar\mF_p$. 
Let $A$ be an abelian variety over $K_0$ and let $X\hookrightarrow A$ be a closed integral subscheme. We shall write $+$ for the group law on $A$. 
  
Let $\Gamma\subseteq A(K_0)$ be a finitely generated subgroup. 

\begin{theor}
Suppose that for any field extension $L_0|K_0$ and any 
$Q\in A(L_0)$, the set \mbox{$X^{+Q}_{L_0}\cap\Tor(A(L_0))$} is not Zariski dense in $X^{+Q}_{L_0}$. Then 
$X\cap\Gamma$ is not Zariski dense in $X$. 
\label{MMprop}
\end{theor}

Here $X^{+Q}_{L_0}$ stands for the scheme-theoretic image of $X_{L_0}$ under 
the morphism $+Q:A_{L_0}\to A_{L_0}$. 
 
\begin{cor} 
Suppose that $X\cap\Gamma$ is Zariski dense in $X$. Then the conclusion 
of the Mordell-Lang conjecture \ref{MLtheor} holds for $F=\bar K_0$, $B=A_{\bar K_0}$ and $Y=X_{\bar K_0}$. 
\label{corimp}
\end{cor}

In an upcoming article, which builds on the present one, C. Corpet (see \cite{Corpet1}) shows that 
Theorem \ref{MMprop} (and thus its corollary) can be generalized; more specifically, he shows that the hypothesis 
that $K_0$ is of transcendence degree $1$ can be dropped, that the hypothesis that 
$\Gamma$ is finitely generated can be weakened to the hypothesis that 
$\Gamma\otimes_\mZ\mZ_{(p)}$ is a finitely generated $\mZ_{(p)}$-module and finally that 
is can be assumed that $A$ is only semiabelian. In particular, he gives a new proof of Theorem \ref{MLtheor}. 

In the present article, we deliberately focus on the situation of an abelian variety and a finitely generated group (which is probably the most important situation) in order to avoid some technical issues, which we feel would obscure the structure of the proof.

The structure of the article is the following. The second section contains 
some general results on the geometry of relative jet schemes (or spaces), which are 
probably known to many specialists but for which there doesn't seem to be a coherent set of references in the literature. The jet spaces considered in \cite{Moosa-Scanlon-Jet} do not seem to suffice 
for our purposes, because they are defined in an absolute situation and the jet spaces 
considered in \cite{Buium-Intersections} are only defined in characteristic $0$ (although this is probably 
not an essential restriction); furthermore, the latter are defined in Buium's language of differential 
schemes, whereas our definition has the philological advantage of being based on the older 
notion of Weil restriction. The subsection 2.1 contains the definition 
of jet schemes and a description of the various torsor structures on the latter. 
The subsection 2.2 contains a short discussion on the structure of the jet schemes of smooth commutative group schemes and various natural maps that are associated with them. In the third section, we use jet schemes to construct some natural schemes in the geometrical context of the  Mordell-Lang conjecture. These "critical 
schemes" are devised to "catch rational points"; we then proceed to show that these  
schemes must be of small dimension. This is deduced from 
a general result on the sparsity of points over finite fields, which are liftable to highly $p$-divisible unramified points. 
This last result is proved in the fourth section. Once we know that the critical schemes are small, it is but a small step to the proof 
of Theorem \ref{MMprop} and Corollary \ref{corimp}. The terminology of the introduction is used in the first three sections but the fourth section has its own terminology and is also technically independent of the rest of the text. A reader who would only be interested in its main result (i.e. Theorem \ref{TVth}) can skip to the fourth section directly. 

The use that we make of jet schemes in this note is in many ways similar 
to the use that A. Buium makes of them in his article on 
the geometric Mordell-Lang conjecture in characteristic $0$ (cf. \cite{Buium-Intersections}). 
In the article \cite{Buium-Voloch-Lang}, where some of Buium's techniques are 
adapted to the context of positive characteristic, the authors give a proof of the Mordell conjecture for curves over 
function fields in positive characteristic, which has exactly the same structure
 as ours, if one leaves out the proof of the result on the the sparsity of liftable points mentioned above. 
 
For more detailed explanations on this connection, see Remarks \ref{CRa} and \ref{CRb} at the end of the text.

{\bf Acknowledgments.} As many people, I am very much indebted to O. Gabber, who pointed out a flaw in an earlier version of this article and who also suggested a way around it. Many thanks to R. Pink for many interesting 
exchanges on the matter of this article. I also want to thank M. Raynaud 
for his reaction on an earlier version of the text and A. Pillay 
for interesting discussions and for suggesting that the method used in this article should work. Finally, I would also like to thank E. Bouscaren and F. Benoist 
for their interest and A. Buium for his very interesting observations. 

\section{Preliminaries.}

We first recall the definition and existence theorem for the Weil restriction functor. Let $T$ be a scheme 
and let $T'\to T$ be a morphism. Let $Z$ be a scheme over 
$T'$. The Weil restriction $\WR_{T'/T}(Z)$ (if it exists) is a $T$-scheme, which represents the functor 
\mbox{$W/T\mapsto\Hom_{T'}(W\times_T T',Z)$.} It is shown in \cite[Par. 7.6]{Bosch-Raynaud-Neron} that $\WR_{T'/T}(Z)$ exists if $T'$ is finite, flat and locally of 
finite presentation over $T$. The Weil restriction is naturally functorial in $Z$ and sends closed immersions to 
closed immersions. The same permanence property is satisfied for smooth and \'etale morphisms. 
Finally notice that the definition of the Weil restriction implies that there is a natural isomorphism 
$\WR_{T'/T}(Z)_{T_1}\simeq \WR_{T'_{T_1}/T_1}(Z_{T'_{T_1}})$ for any scheme $T_1$ over $T$ (in words: Weil restriction 
is invariant under base-change on $T$). See \cite[Par. 7.6]{Bosch-Raynaud-Neron}  for all this. 

\subsection{Jet schemes}
 
Let $k_0$ be field. 
Let $U$ be a smooth scheme over $k_0$. Let $\Delta:U\to U\times_{k_0} U$ be the diagonal immersion. 
Let $I_\Delta\subseteq\CO_{U\times_{k_0} U}$ be the ideal sheaf of $\Delta_*U$. 
For all $n\in\mN$, we let $U_n:=\CO_{U\times_{k_0} U}/I_\Delta^{n+1}$ be the $n$-th infinitesimal neighborhood of 
the diagonal in $U\times_{k_0} U$. 

Write $\pi_1,\pi_2:U\times_{{k_0}} U\to U$ for the first and second projection 
morphism, respectively. Write $\pi_1^{U_n},\pi_2^{U_n}:U_n\to U$ for the induced morphisms. 
We view $U_n$ as a $U$-scheme via the {\it first} projection $\pi_1^{U_n}$. 

We write $i_{m,n}:U_m\hookrightarrow U_n$ for the natural inclusion morphism. 

\begin{lemma}
The $U$-scheme $U_n$ is flat and finite.
\end{lemma}
\beginProof
As a $U$-scheme, $U_n$ is finite because it is quasi-finite and proper over $U$, since 
$U^\red_n=\Delta_*(U)$. So we only have to prove that 
it is flat over $U$. For this purpose, we may view $U_n$ as a coherent sheaf of $\CO_U$-algebras 
(via the second projection). 

Let $I:=I_\Delta$. For any $n\geqslant 0$, there are exact sequences of $\CO_{U_{n+1}}$-modules 
(and hence $\CO_U$-modules) 
$$
0\to I^{n+1}/I^{n+2}\to\CO_{U_{n+1}}\to\CO_{U_{n}}\to 0
$$
Furthermore $I^{n+1}/I^{n+2}$ is naturally a $\CO_{U_0}$-module and is isomorphic to $\Sym^{n+1}_{\CO_{U_1}}(I/I^2)$ as a $\CO_{U_0}$-module, 
because $I$ is locally generated by a regular sequence 
in $U\times_{{k_0}}U$ ($U$ being smooth over ${k_0}$). See 
\cite[chap. 6, 16.]{Matsumura-Commutative} for this. Hence $I^{n+1}/I^{n+2}$ is locally free as 
a $\CO_U$-module. Since $U_0=\Delta_*(U)$ is locally free as 
a $\CO_U$-module, we may apply induction on $n$ to prove that $\CO_{U_{n}}$ is locally free, which 
is the claim.\endProof

Let $W/U$ be a scheme over $U$. 

\begin{defin}
The $n$-th jet scheme $J^n(W/U)$ of $W$ over $U$ is the $U$-scheme 
$\WR_{U_n/U}(\pi_2^{U_n,*}W)$.
\end{defin}
By $\pi_2^{U_n,*}W$ we mean the base-change of $W$ to $U_n$ via the 
morphism $\pi_2^{U_n}:U_n\to U$ described above. 

If $W_1$ is another scheme over $U$ and $W\to W_1$ is a morphism of $U$-schemes, then the induced morphism 
$\pi_2^{U_n,*}W\to\pi_2^{U_n,*}W_1$ over $U_n$ leads to 
a morphism of jet schemes $J^n(W/U)\to J^n(W_1/U)$ over $U$, so that 
the construction of jet schemes is covariantly functorial in $W$. 

Notice that the permanence properties of Weil restrictions show that if the morphism 
$W\to W_1$ 
is a closed immersion, then so is the morphism $J^n(W/U)\to J^n(W_1/U)$. 
Same for smooth and \'etale morphisms. 

To understand the nature of jet schemes better, let $u\in U$ be a closed point. Suppose until the end 
of the sentence following the string of equations \refeq{jetpoint} that $k_0$ is algebraically closed. View $u$ as closed 
reduced subscheme of $U$. Let $u_n$ be the $n$-th infinitesimal neighborhood of $u$ in 
$U$. From 
the definitions, we infer that there are canonical bijections
\begin{eqnarray}
J^n(W/U)(u)&=&J^n(W/U)_u({k_0})=\Hom_{U_n}(u\times_U U_n,\pi_2^{U_n,*}W)\nonumber\\
&=&\Hom_{U_n}(u_n,W_{u_n})=\Hom_{u_n}(u_n,W_{u_n})=W(u_n)
\label{jetpoint}
\end{eqnarray}
In words, \refeq{jetpoint} says the set of geometric points of the fibre of $J^n(W/U)$ over $u$ 
corresponds to the set of sections of $W$ over the $n$-th infinitesimal neighborhood of $u$; the scheme 
$J^n(W/U)_u$ is often called the scheme of arcs of order $n$ at $u$ in the literature (see \cite[Ex. 2.5]{Moosa-Scanlon-Jet}). 

The family of $U$-morphisms $i_{m,n}:U_m\to U_{n}$ induce 
$U$-morphisms \mbox{$\Lambda^W_{n,m}:J^{n}(W/U)\to J^{m}(W/U)$} for any \mbox{$m\leqslant n$.} These morphisms will be studied in detail in the proof of the next lemma. 

\begin{lemma} 
Suppose that $W$ is a smooth $U$-scheme. 
For all $n\geqslant 1$, 
the morphism $$\WR_{U_n/U}(\pi_2^{U_n,*}W)\to \WR_{U_{n-1}/U}(\pi_2^{U_{n-1},*}W)$$ 
makes $\WR_{U_n/U}(\pi_2^{U_n,*}W)$ into 
a $\WR_{U_{n-1}/U}(\pi_2^{U_{n-1},*}W)$-torsor under the vector bundle\linebreak\mbox{$\Lambda^{W,*}_{n,0}(\Omega_{W/U}^\vee)
\otimes\Sym^{n}({\Omega}_{U/{k_0}})$.}
\label{torslem}
\end{lemma}
\beginProof  
Let $T\to U$ be an affine $U$-scheme. By definition
$$
\WR_{U_n/U}(\pi_2^{U_{n},*}W)(T)\simeq \Hom_{U_n}(T\times_{U}
U_n,\pi_2^{U_{n},*}W)
$$
and
$$
\WR_{U_{n-1}/U}(\pi_2^{U_{n-1},*}W)(T)\simeq \Hom_{U_{n-1}}(T\times_{U}U_{n-1},\pi_2^{U_{n-1},*}W).
$$
Now the immersion $U_{n-1}\hookrightarrow U_n$ gives rise to a natural restriction map 
\begin{equation}
\Hom_{U_n}(T\times_{U}U_n,\pi_2^{U_{n},*}W)\to \Hom_{U_{n-1}}(T\times_{U}U_{n-1},\pi_2^{U_{n-1},*}W).
\label{urss}
\end{equation}
This is the functorial description of the morphism $\WR_{U_n/U}(\pi_2^{U_{n},*}W)\to \WR_{U_{n-1}/U}(\pi_2^{U_{n-1},*}W)$. 

Now notice that the ideal of  $U_{n-1}$ in $U_n$ is a square $0$ ideal. 

Let $f\in \Hom_{U_{n-1}}(T\times_{U}U_{n-1},\pi_2^{U_{n-1},*}W)$. View $f$ as a $U_n$-morphism 
$T\times_{U}U_{n-1}\to \pi_2^{U_n,*}W$ via the canonical closed immersions 
$\pi_2^{U_{n-1},*}W\hookrightarrow \pi_2^{U_n,*}W$ and $U_{n-1}\hookrightarrow U_n$. 
The fibre over $f$ of the map \refeq{urss} then consists of the extensions of $f$ to $U_n$-morphisms 
\mbox{$T\times_{U}U_n\to\pi_2^{U_n,*}W$.} 
The theory of infinitesimal extensions of morphisms to smooth schemes (see \cite[Exp. III, Prop. 5.1]{SGA1})
 implies that this fibre 
is an affine space under the group 
$$
H^0(T\times_{U}U_{n-1}, f^*\Omega_{\pi_2^{U_{n},*}W/U_n}^\vee\otimes N)
$$
where $N$ is the conormal bundle of the closed immersion $T\times_{U}U_{n-1}\hookrightarrow 
T\times_{U}U_n$. Since $U_n$ and $U_{n-1}$ are flat over $U$, the coherent sheaf 
$N$ is the pull-back to $T\times_{U}U_{n-1}$ of the conormal bundle of the immersion 
$U_{n-1}\hookrightarrow U_n$. Now since the diagonal is regularly immersed in $U\times_{{k_0}} U$ (because $U$ is smooth over ${k_0}$), 
the conormal bundle of the immersion ${U}_{n-1}\hookrightarrow {U}_{n}$ is 
$\Sym^{n}({\Omega}_{U/{k_0}})$ (viewed as 
a sheaf in $\CO_{U_{n-1}}$-modules via the closed immersion $U_0\to U_{n-1}$). See 
\cite[chap. 6, 16.]{Matsumura-Commutative} for this. Hence 
\begin{eqnarray*}
&&H^0(T\times_{U}U_{n-1}, f^*\Omega_{\pi_2^{U_{n},*}W/U_n}^\vee\otimes N)\\
&\simeq& 
H^0(T\times_{U}U_{n-1}, f^*\Omega_{\pi_2^{U_{n},*}W/U_n}^\vee\otimes
\Sym^{n}({\Omega}_{U/{k_0}}))
\simeq H^0(T, f^*_0\Omega_{W/U}^\vee\otimes\Sym^{n}({\Omega}_{U/{k_0}}))
\end{eqnarray*}
where $f_0$ is the $U$-morphism $T\to W$ arising from $f$ by base-change to $U$. 
This proves the lemma.\endProof

\subsection{The jet schemes of smooth  commutative group schemes}
 
We keep the terminology of the last subsection. 
Let $\CC/U$ be a commutative group scheme over $U$, with 
zero-section $\epsilon:U\to\CC$. If $n\in\mN$, we shall write 
$[n]_\CC:\CC\to\CC$ for the multiplication-by-$n$ morphism. The schemes  
$J^n(\CC/U)$ are then 
naturally group schemes over $U$.  
Furthermore, for each $n\geqslant m\geqslant 0$, the morphism $\Lambda^\CC_{n,m}:J^{n}(\CC/U)\to J^{m}(\CC/U)$ is a 
 morphism of group 
schemes. If $m=n-1$, the kernel of $\Lambda^\CC_{n,m}$ is the vector bundle \mbox{$\epsilon^*(\Omega_{\CC/U}^\vee)
\otimes \Sym^{n}({\Omega}_{U/{k_0}})$.} The torsor structure 
is realized by the natural action of \mbox{$\epsilon^*(\Omega_{\CC/U}^\vee)
\otimes \Sym^{n}({\Omega}_{U/{k_0}})$} on  $J^n(\CC/U)$. 
The details of the verification of these facts are left to the reader. 

\begin{lemma} Let $n\geqslant 1$. Suppose that $\char(k_0)=p$. There is an $U$-morphism 
$"p^{n}":\CC\to J^{n}(\CC/U)$ such that \mbox{$\Lambda^\CC_{n,0}\circ"p^{n}"=[p^{n}]_\CC$} and 
$[p^{n}]_{J^n(\CC/U)}="p^{n}"\circ\Lambda^\CC_{n,0}$.
\label{pquotelem}
\end{lemma}
\beginProof
Let $T\to U$ be an affine $U$-scheme. Define a map  
$$
\phi_{T,n}:\Hom_U(T,\CC)\to\Hom_{U_n}(T\times_U U_n,\pi_2^*\CC)
$$
by the following recipe. Let $f\in \Hom_U(T,\CC)$ and take any extension 
$\widetilde{f}$ of $f$ to a morphism 
$T\times_U U_n\to (\pi_2^*\CC)_{U_n}$; 
then define $\phi_{T,n}(f)=p^{n}\cdot\widetilde{f}$. To see that this does not depend on 
the choice of the extension $\widetilde{f}$, notice that the kernel $K_n$ of the restriction map 
$$
\Hom_{U_n}(T\times_U U_n,\pi_2^*\CC)\to \Hom_U(T,\CC)
$$
is obtained by successive extensions by the groups 
$
H^0(T, f^*\Omega_{\CC/U}^\vee\otimes\Sym^{i}({\Omega}_{U/{k_0}})
$ for $i=1,\dots,n$ (see \cite[Chap. III, 5., Cor. 5.3]{SGA1} for all this). Hence $K_n$ is annihilated by multiplication by 
$p^{n}$ because $T$ is a scheme of characteristic $p$. 

The definition of $\phi_{T,n}$ is functorial in $T$ and thus by patching 
the morphisms $\phi_{T,n}$ as $T$ runs over the elements of an affine cover of 
$\CC$ we 
obtain the required morphism $"p^{n}"$.
\endProof

Now notice that for any scheme $W$ over $U$, there is a canonical map $\lambda_n^W:W(U)\to J^{n}(W/U)(U)$, which sends the $U$-morphism 
\mbox{$f:U\to W$} to $J^n(f):J^n(U/U)=U\to J^n(W/U)$. 

 \begin{lemma}
The maps $\lambda^W_n$ have the following properties : 
\begin{itemize}
\item[\rm (a)] for $n\geqslant m\geqslant 0$ the identity $\Lambda^W_{n,m}\circ\lambda^W_n=\lambda^W_m$; \item[\rm (b)] if $W/U$ is commutative group scheme over $U$, then $\lambda^W_n$ is 
a homomorphism; furthermore on $W(U)$ we then have the identity \mbox{$[p^{n}]_{J^n(W/U)}\circ\lambda_n^W="p^{n}"$}; 
\item[\rm (c)] if $f:W\to W_1$ is a $U$-scheme morphism then $J^n(f)\circ\lambda^W_n=
\lambda^{W_1}_n\circ f$.
\end{itemize}
\label{cruclem}
\end{lemma}
\beginProof (of Lemma \ref{cruclem}). Exercise for the reader. \endProof

\begin{remark}\textnormal{   An interesting feature of the map $\lambda_n^W$ is that it does \underline{not} arise 
from a morphism of schemes $W\to J^{n}(W/U)$. }
\end{remark}

\section{Proof of Theorem \ref{MMprop} \& Corollary \ref{corimp}}

We now turn to the proof of our main result. We shall use the terminology of the preliminaries.
Let $k_0:=\bar\mF_p$ and suppose now that $U$ is a smooth curve over $k_0$, whose 
function field is $K_0$. 
We take $U$ sufficiently small so that $X$ extends to a flat scheme $\CX$ over $U$ and 
so that $A$ extends to an abelian scheme $\CA$ over $U$. We also suppose that 
the closed immersion $X\hookrightarrow A$ extends to a closed immersion $\CX\to\CA$. 

Recall that the following hypothesis is supposed to hold : for any field extension $L_0|K_0$ and any 
$Q\in A(L_0)$, the set $X^{+Q}_{L_0}\cap\Tor(A(L_0))$ is not Zariski dense in $X^{+Q}_{L_0}.$

\subsection{The critical schemes}

For all $n\geqslant 0$, we define 
$$\Crit^n(\CX,\CA):=[p^{n}]_*(J^{n}(\CA/U))\cap J^{n}(\CX/U).$$ 

Here $[p^{n}]_*(J^{n}(\CA/U))$ is the scheme-theoretic image of $J^{n}(\CA/U)$ by $[p^{n}]_{J^{n}(\CA/U)}.$ 
Notice that by Lemma \ref{pquotelem}, we have $[p^{n}](J^{n}(\CA/U))="p^{n}"(\CA)$ and since $[p^{n}]$ is proper (because $\CA$ is proper over $U$), we see that $[p^{n}](J^{n}(\CA/U))$ is 
closed and that the natural morphism \mbox{$[p^{n}]_*(J^{n}(\CA/U))\to\CA$} is finite. 

The morphisms $\Lambda^\CA_{n,n-1}:J^n(\CA/U)\to J^{n-1}(\CA/U)$ lead to a projective system of 
$U$-schemes
$$
\dots\to\Crit^{2}(\CX,\CA)\to\Crit^1(\CX,\CA)\to\CX.
$$
whose connecting morphisms are finite. We let $\Exc^n(\CA,\CX)\hookrightarrow\CX$ be the 
scheme-theoretic image of $\Crit^n(\CA,\CX)$ in $\CX$. 


For any $Q\in\CA(U)=A(K_0)$, we shall write $\CX^{+Q}=\CX+Q$ for the translation of $\CX$ by $Q$ in $\CA$. 
\begin{prop}
There exists $\alpha=\alpha(\CA,\CX)\in\mN$ such that for all $Q\in\Gamma$ the set $\Exc^\alpha(\CA,\CX^{+Q})$ is not 
dense in $\CX^{+Q}$. 
\label{critprop}
\end{prop}

\begin{remark}\textnormal{Proposition \ref{critprop} should be compared 
to Theorem 1 in \cite{Buium-Intersections}.}\end{remark}

The following theorem, proved by galois-theoretic methods in section \ref{sec4}, will play a crucial role in the proof of Proposition \ref{critprop}. 

Let $S:=\Spec\, k_0[[t]]$. Let $L:=k_0((t))$ be the function field of $S$. For any $n\in\mN$, let 
$S_n:=\Spec\,k_0[t]/t^{n+1}$ be the $n$-th infinitesimal neighborhood of the closed point of $S$ in $S$. Fix $\lambda_0\in\mN^*$ and let 
$R^\alg=R^{\alg,\lambda_0}:=\mF_{p^{\lambda_0}}[[t]]\subseteq  k_0[[t]].$ Let $S^\alg=S^{\alg,\lambda_0}:=\Spec\,R^\alg.$ There is an obvious morphism 
$S\to S^\alg$. 

Let $\CD$ be an abelian scheme over $S$ and let $\CZ\hookrightarrow\CD$ be a closed integral subscheme. 
Suppose that the abelian scheme has a model $\CD^\alg$ over $S^\alg$ as an abelian scheme and that 
the immersion  $\CZ\hookrightarrow\CD$ has a model $\CZ^\alg\hookrightarrow\CD^\alg$ over $S^\alg.$ 
If $c\in \CD(S)$, write as usual $\CZ^{+c}:=\CZ+c$ for the translation of 
$\CZ$ by $c$ in $\CD$. Let $D_0$ (resp. $D$) be the fibre of $\CD$ over the closed point of $S$ (resp. over the generic point of $S$). 
If $c\in\CD(S)$, let $Z_0^{+c}$ (resp. $Z^{+c}$) be the fibre of $\CZ^{+c}$ over the closed point of $S$ (resp. over the generic point of $S$). 

Notice that there is a natural inclusion $\CD^\alg(S^\alg)\subseteq \CD(S)$. 

\begin{theor}
Suppose that $\Tor(D(\bar L))\cap X^{+c}_{\bar L}$ is not dense in $X^{+c}_{\bar L}$ for all \mbox{$c\in \CD^\alg(S^\alg)\subseteq \CD(S)$.}

Then there exists a constant $n_0=n_0(\CD,\CZ)\in\mN^*$, such that for all $c\in\CD^\alg(S^\alg)\subseteq \CD(S)$, the set 
$$
\{P\in Z_0^{+c}(k_0))\ |\ \textnormal{\rm $P$ lifts to an element of $\CZ^{+c}(S_{n_0})\cap p^{n_0}\cdot
\CD(S_{n_0})$}\}
$$
is not Zariski dense in $Z_0^{+c}$. 
\label{pffth}
\end{theor}
\beginProof (of Theorem \ref{pffth}). This is a special case of Corollary \ref{TVuni}.\endProof

\beginProof (of Proposition \ref{critprop}). 
Since $\CX$ is flat over $U$ and $X$ is integral, we see that 
$\CX$ is also integral (see for instance \cite[4.3.1, Prop. 3.8]{Liu-Algebraic} for this). Hence it is sufficient 
to show that $\Exc^n(\CA,\CX^{+Q})_u$ is not Zariski dense in $\CX_u^{+Q}$ for some (any) closed point 
$u\in U$. Now using \refeq{jetpoint} in the previous section, we see that 
\begin{eqnarray*}
\Crit^n(\CA,\CX^{+Q})_u({k_0})&=&
([p^{n}]_*(J^{n}(\CA/U)))_u({k_0})\cap J^{n}(\CX^{+Q}/U)_u({k_0})\\
&=&\{P\in J^{n}(\CX^{+Q}/U)_u({k_0})\,|\, \exists \widetilde{P}\in J^{n}(\CA/U)_u({k_0})\,:\,p^{n}\cdot\widetilde{P}=P\}\\
&=&\{P\in \CX^{+Q}(u_n)\,|\, \exists \widetilde{P}\in \CA(u_n)\,:\,p^{n}\cdot\widetilde{P}=P\}
\end{eqnarray*}
and thus
\begin{eqnarray*}
\Exc^n(\CA,\CX^{+Q})_u({k_0})&=&
\{P\in\CX_u^{+Q}({k_0})\,|\,\textrm{$P$ lifts to an element of $\CX^{+Q}(u_n)\cap p^{n}\cdot\CA(u_n)$}\}
\end{eqnarray*}
Now notice that $\CA$ has a model $\widetilde{\CA}$ as an abelian scheme over a curve $\widetilde{U}$, which is smooth over 
a finite field; also since the group $\Gamma$ is finitely 
generated, we might assume that $\Gamma$ is the image of a group $\widetilde{\Gamma}\subseteq\widetilde{\CA}(\widetilde{U})$. Finally, 
we might assume that the immersion $\CX\hookrightarrow\CA$ has a model $\widetilde{\CX}\hookrightarrow\widetilde{\CA}$ over $\widetilde{U}$. 
We may thus apply Theorem \ref{pffth} to the base-change of $\CX\hookrightarrow\CA$ to the completion of $U$ at $u$. 
We obtain  that there is an $n_0$ such that the set 
$$
\{P\in\CX_u^{+Q}({k_0})\,|\,\textrm{$P$ lifts to an element of $\CX^{+Q}(u_{n_0})\cap p^{n_0}\cdot\CA(u_{n_0})$}\}
$$
is not Zariski dense in $\CX_u$ for all $Q\in\Gamma$. So we may 
set $\alpha=n_0$.\endProof

\subsection{End of proof}

The proof of Theorem \ref{MMprop} is by contradiction. So suppose that $X\cap\Gamma$ is dense in $X$. 

Let $P_1\in \Gamma$ be such that $(X+P_1)\cap p\cdot \Gamma$ is dense, let $P_2\in p\cdot \Gamma$ such that 
$(X+P_1+P_2)\cap p^2\cdot \Gamma$ is dense in $X+P_1+P_2$ and so forth. The existence of 
the sequence of point $(P_i)_{i\in\mN^*}$ is guaranteed by the assumption on $\Gamma$, which 
implies that $p^i\Gamma/p^{i+1}\Gamma$ is finite for all $i\geqslant 0$. 

Now let $\alpha=\alpha(\CA,\CX)$ be the natural number provided by Proposition \ref{critprop}. 
Let $Q=\sum_{i=1}^{\alpha}P_i$. By construction, the set 
$\CX^{+Q}\cap p^{\alpha}\cdot \Gamma$ is dense in $\CX^{+Q}$. On the other hand, by Lemma \ref{cruclem}, 
\begin{eqnarray*}
\CX^{+Q}(U)\cap p^{\alpha}\cdot \Gamma&=&
\Lambda^\CA_{\alpha,0}(\lambda^\CA_\alpha(\CX^{+Q}(U)\cap p^\alpha\cdot \Gamma)) 
\subseteq \Lambda^\CA_{\alpha,0}[\,\lambda^\CX_\alpha(\CX^{+Q}(U))\cap\lambda^\CA_\alpha(p^\alpha\cdot \Gamma)\,]\\
&\subseteq& \Lambda^\CA_{\alpha,0}[\,J^\alpha(\CX^{+Q}/U)\cap p^{\alpha}\cdot J^\alpha(\CA/U)(U)\,]
\subseteq\Lambda^\CX_{\alpha,0}[\,\Crit^\alpha(\CA,\CX^{+Q})\,]\\
&=&\Exc^\alpha(\CA,\CX^{+Q})
\end{eqnarray*}
and thus we deduce that $\Exc^\alpha(\CA,\CX^{+Q})$ is dense in $\CX^{+Q}$. This contradicts 
Proposition \ref{critprop} and concludes the proof of Theorem \ref{MMprop}.

The proof of Corollary \ref{corimp} now follows directly from Theorem \ref{MMprop} and from the following invariance lemma.  

\begin{lemma}
Suppose that the hypotheses of Theorem \ref{MLtheor} hold.
Let $F'$ be an algebraically closed field and let $F'|F$ be a field extension. 
 Then Theorem \ref{MLtheor} holds if and only if Theorem 
 \ref{MLtheor} holds, with $F'$ in place of $F$, $Y_{F'}\hookrightarrow B_{F'}$ in place 
 of  $Y\hookrightarrow B$, and the image $\Lambda_{F'}\subseteq B_{F'}(F')$ of $\Lambda$ in place of 
 $\Lambda$. 
 \label{invlem}
 \end{lemma}
 \beginProof
 The implication $\Longrightarrow$ follows from the fact that $Y_{F'}\cap\Lambda_{F'}$ is dense 
 in $Y_{F'}$ if and only if  $Y\cap\Lambda$ is dense 
 in $Y$; indeed the morphism  $\Spec\,F'\to\Spec\,F$ is universally open (see \cite[IV, 2.4.10]{EGA} for this).
  
 Now we prove the implication $\Longleftarrow$. Let $C_1:=\Stab(Y_{F'})^\red$ and suppose that 
 there exists 
 \begin{itemize}
\item a semiabelian variety $B'_1$ over $F'$; 
\item a homomorphism with finite kernel $h_1: B'_1\to B_{F'}/C_1$;
\item a model ${\bf B}'_1$ of $B'_1$ over a finite subfield $\mF_{p^r} \subset F'$;
\item 
an irreducible reduced closed subscheme ${\bf Y}'_1\hookrightarrow{\bf B'}_1$;
\item a point $b_1\in (B_{F'}/C_1)(F')$, such that\,\, 
$Y_{F'}/C_1 = b_1 + h_{1,*}\bigl({\bf Y}'_1\times_{\mF_{p^r}}F'\bigr).$
\end{itemize}
Now first notice that since $\Stab(\bullet)$ represents a functor, there is a natural isomorphism 
$\Stab(Y_{F'})\simeq\Stab(Y)_{F'}$ and since $F$ is algebraically closed also a natural isomorphism 
$\Stab(Y_{F'})^\red\simeq(\Stab(Y)^\red)_{F'}$. Secondly, we have  $\mF_{p^r} \subset F$, since 
$F$ is algebraically closed. Thirdly, if $B_2$ and $B_3$ are semiabelian varieties over $F$ and 
$\phi:B_{2,F'}\to B_{3,F'}$ is a homomorphism of group schemes over $F'$, then $\phi$ arises by base-change from 
an $F$-morphism $B_2\to B_3$. This is a consequence of the fact that the graph of $\phi$ has a dense set of 
torsion points in $B_{2,F'}\times_{F'}B_{3,F'}$, and torsion points are defined in $B_{2}\times_{F}B_{3}$. 
Putting these facts together, we deduce that there exists
\begin{itemize}
\item a semiabelian variety $B'$ over $F$; 
\item a homomorphism with finite kernel $h: B'\to B/C$;
\item a model ${\bf B}'$ of $B'$ over a finite subfield $\mF_{p^r} \subset F$;
\item 
an irreducible reduced closed subscheme ${\bf Y}'\hookrightarrow{\bf B}'$;
\item a point $b_1\in (B_{F'}/C_{F'})(F')$, such that\,\, 
$Y_{F'}/C_{F'} = b_1 + h_{F',*}\bigl({\bf Y}'\times_{\mF_{p^r}}F'\bigr).$
\end{itemize}
where $C=\Stab(Y)^\red$. Now last point in the list above shows that 
\mbox{$\Transp(Y_{F'}/C_{F'}, h_{F',*}\bigl({\bf Y}'\times_{\mF_{p^r}}F'\bigr))(F')\not=\emptyset$.} Here 
$\Transp(\bullet)$ is the transporter, which is a generalization of the stabilizer (see \cite[Exp. VIII, 6.]{SGA3-2} for the 
definition). Thus
 $\Transp(Y/C, h_{*}\bigl({\bf Y}'\times_{\mF_{p^r}}F\bigr))(F)\not=\emptyset$, which is to say that there also exists 
 \begin{itemize}
\item a point $b_1\in (B/C)(F)$, such that\,\, 
$Y/C = b_1 + h_{*}\bigl({\bf Y}'\times_{\mF_{p^r}}F\bigr).$
\end{itemize}
 This concludes the proof.\endProof

\section{Sparsity of highly $p$-divisible unramified liftings}

\label{sec4}

This section can be read independently of the rest of the text and its results do not rely on the 
previous ones. Also, unlike the previous sections, {\it the terminology of this section is independent of the terminology of the introduction.}  

Let $S$ be the spectrum of complete discrete local ring. Let  $k$ be the residue field of its closed point. We suppose that 
$k$ is a {\it finite field} of characteristic $p$. 
Let $K$ be the fraction field of $S$. Let $S^\sh$ be the spectrum of the strict henselisation of $S$ and let 
$L$ be the fraction field of $S^\sh$. We identify $\bar k$ with the residue field of the closed point of $S^\sh$. 
For any $n\in\mN$, we shall write $S_n$ (resp. $S^\sh_n$) for the $n$-th infinitesimal neighborhood of the closed point 
of $S$ (resp. $S^\sh$) in $S$ (resp. $S^\sh$). 

Let $\CA$ be an abelian scheme over $S$  and let $A:=\CA_K$. Write $A_0$ for the fibre of $\CA$ over the closed 
point of $S$. 

\begin{theor}
Let $\CX\hookrightarrow\CA$ be a closed integral subscheme. 
Let $X_0$ be the fibre of $\CX$ over the closed point of $S$ and let $X:=\CX_K$. 

Suppose that $\Tor(A(\bar K))\cap X_{\bar K}$ is not dense in $X_{\bar K}$.

Then there exists a constant $m\in\mN$, such that the set 
$$
\{P\in X_0(\bar k))\ |\ \textrm{\rm $P$ lifts to an element of $\CX(S^\sh_{m})\cap p^{m}\cdot
\CA(S^\sh_{m})$}\}
$$
is not Zariski dense in $X_0$.
\label{TVth}
\end{theor}

Suppose for the time of the next sentence that $S$ is the spectrum of a complete 
discrete valuation ring, which is absolutely unramified and is the 
completion of a number field along a non-archimedean place. In this situation, M. Raynaud proves Theorem 
\ref{TVth} and Corollary \ref{TVuni} below, under the stronger hypothesis that 
$X_{\bar K}$ does not contain any translates of positive-dimensional abelian subvarieties 
of $A_{\bar K}$ (see \cite[Prop. II.1.1]{Raynaud-Around}). See also \cite[Th. II, p. 207]{Raynaud-Courbes} for 
a more precise result in the situation where $X$ is a smooth curve. 

In the case where $S$ is the spectrum of the ring of integers of a finite extension of 
$\mQ_p$, Theorem \ref{TVth} implies versions of the Tate-Voloch conjecture 
(see \cite{Tate-Voloch-Linear} and \cite{Scanlon-The-conjecture}). We leave it to the reader to work out the details. 

Preliminary to the proof of \ref{TVth}, we quote the following result. 
Let $B$ be an abelian variety over an algebraically closed field $F$ and let $\psi:B\to B$ be an endomorphism. Let $R\in\mZ[T]$ be a polynomial, which has no roots of unity among its complex roots. 
Suppose that $R(\psi)=0$ in the ring of endomorphisms of $B$. 

\begin{prop}[Pink-R\"ossler]
Let $Z\subseteq B$ be a closed irreducible subset such that $\psi(Z)=Z$. Then 
$\Tor(B(F))\cap Z$ is dense in $Z$.
\label{prop61PR2}
\end{prop}

The proof of Proposition \ref{prop61PR2} is based on a spreading out argument, which is used to reduce the problem 
to the case where $F$ is the algebraic closure of a finite field.  In this last case, the statement becomes obvious. 
See \cite[Prop. 6.1]{PR2} for the details. 

We shall use the map $"p^l":A_0(\bar k)\to\CA(S^\sh_{l})$, which is defined 
by the formula $"p^l"(x)=p^l\cdot\widetilde{x}$, where $\widetilde{x}$ is any lifting of $x$ (this does not depend 
on the lifting; see \cite[after Th. 2.1]{Katz-Serre-Tate}).

\beginProof (of Theorem \ref{TVth}).   Let $\phi$ be a topological 
generator of $\Gal(\bar k|k)$. By the Weil conjectures for abelian varieties, there is a a polynomial 
$$Q(T):=T^{2g}-(a_{2g-1} T^{2g-1}+\dots+ a_0)$$ with $a_i\in\mZ$, such that 
$Q(\phi)(x)=0$ for all $x\in A_0(\bar k)$ and such that $Q(T)$ has no roots of unity among its complex roots. 
Let $M$ be the matrix 
\medskip
$$
\begin{bmatrix}
0&1&\ldots&0&0\cr
\vdots&\vdots&&\vdots&\vdots\cr
0&0&\ldots&0&1\cr
a_0&a_1&\ldots&a_{n-2}&a_{2g-1}\cr
\end{bmatrix}
$$
\medskip
We view $M$ as an endomorphism of abelian $S$-schemes $\CA^{2g}\to\CA^{2g}$. 
Let \mbox{$\tau\in\Aut_{S}(S^\sh)$} be the canonical lifting of $\phi$. By construction, 
$\tau$ induces an element of $\Aut_{S_n}(S^\sh_n)$ for any $n\geqslant 0$, which we also call $\tau$. 
The reduction map $\CA(S^\sh)\to \CA(S_n^\sh)$ is compatible with the action of $\tau$ on both sides.
 Write $$u(x):=(x,\tau(x),\tau^2(x),\dots,\tau^{2g-1}(x))\in({\prod^{2g-1}_{s=0}}\CA)(S^\sh)$$ for any element 
$x\in \CA(S^\sh)$. Abusing notation, we shall also write
$$u(x):=(x,\tau(x),\tau^2(x),\dots,\tau^{2g-1}(x))\in({\prod^{2g-1}_{s=0}}\CA)(S^\sh_n)$$ 
for any element 
$x\in \CA(S^\sh_n)$. 
 By construction, for any $x\in\CA(S^\sh)$ (resp. any $x\in\CA(S_n^\sh)$), 
the equation $Q(\tau)(x)=0$ implies the vector identity $M(u(x))=u(\tau(x))$. 

Now consider the closed $S$-subscheme of $\CA^{2g}$
$$
\TV:=\cap_{t\geqslant 0}M^{t}_*{\Big(}\cap_{r\geqslant 0}M^{r,*}(\prod^{2g-1}_{s=0}\CX){\Big)}
$$
where for any $r\geqslant 0$, $M^r$ is the $r$-th power of $M$. The symbol $M^{t}_*(\cdot)$ refers to the 
scheme-theoretic image and the intersections are the scheme-theoretic intersections. The intersections are finite by noetherianity. 

Let $\lambda:J\to\CA^{2g}$ be a morphism of schemes. The construction of $\TV$ implies that if 
the following conditions are verified 

(i) $M^r\circ \lambda$ factors through $\prod^{2g-1}_{s=0}\CX$ for all $r\geqslant 0$ and 

(ii) for all $r\geqslant 0$, there is a morphism $\phi_r:J\to\cap_{r\geqslant 0}M^{r,*}(\prod^{2g-1}_{s=0}\CX)$ such that 
$\lambda=M^r\circ\phi_r$ 

then $\lambda$ factors through $\TV$. 

In particular, if (i) is verified and $M^{r_\lambda}\circ\lambda=\lambda$ 
for some $r_\lambda\geqslant 1$, then 
$\lambda$ factors through $\TV$. 

\begin{remark}\textnormal{In particular, this implies the following:  if $x\in\CX(S^\sh)$ (resp. $x\in\CX(S_n^\sh)$) has the property that $Q(\tau)(x)=0$, 
then $u(x)\in\TV(S^\sh)$ (resp. $u(x)\in \TV(S^\sh_n)$).}
\label{remQ}
\end{remark}

\begin{lemma}
There is a set-theoretic identity $M(\TV)=\TV$.
\end{lemma}
\beginProof (of the lemma)
Since $M$ is proper, we have a set-theoretic identity
$$
\TV=\cap_{t\geqslant 0}M^{t}{\Big(}\cap_{r\geqslant 0}M^{r,-1}(\prod^{2g-1}_{s=0}\CX){\Big)}
$$
Now directly from the construction, we have 
$$
M{\Big(}\cap_{r\geqslant 0}M^{r,-1}(\prod^{2g-1}_{s=0}\CX){\Big)}\subseteq\cap_{r\geqslant 0}M^{r,-1}(\prod^{2g-1}_{s=0}\CX)
$$
and hence we have inclusions
$$
\cap_{r\geqslant 0}M^{r,-1}(\prod^{2g-1}_{s=0}\CX)\supseteq M{\Big(}\cap_{r\geqslant 0}M^{r,-1}(\prod^{2g-1}_{s=0}\CX){\Big)}
\supseteq M^2{\Big(}\cap_{r\geqslant 0}M^{r,-1}(\prod^{2g-1}_{s=0}\CX){\Big)}\supseteq\dots 
$$
and thus by noetherianity 
$$
M^{l}{\Big(}\cap_{r\geqslant 0}M^{r,-1}(\prod^{2g-1}_{s=0}\CX){\Big)}
=M^{l+1}{\Big(}\cap_{r\geqslant 0}M^{r,-1}(\prod^{2g-1}_{s=0}\CX){\Big)}
$$
for some $l\geqslant 0$. This implies the result. 
\endProof

Now we apply Proposition \ref{prop61PR2} and we obtain that $\TV_{\bar K,\red}\cap\Tor(\prod^{2g-1}_{s=0}A(\bar K))$ is dense in $\TV_{\bar K,\red}$. Hence the projection onto the first factor 
$\TV_{K}\to X$ is not surjective by hypothesis. 

 Let $T$ be the scheme-theoretic image of the morphism $\TV\to\CX$ given by the first 
 projection. Notice that $X_0$ is a closed subscheme of $T$, because every element $P$ of $X_0(\bar k)$ satisfies the equation 
$Q(\phi)(P)=0$. Let $H$ be the closed subset of $T$, which is the union of the irreducible components of $T$, which surject onto $S$. A reduced irreducible component $I$ of $T$, which surjects onto $S$, is flat over $S$; since $H\not=\CX$, we have in particular $I\not=\CX$ and so we see that 
 the dimension of the fibre of $I$ over the closed point of $S$ is strictly smaller than the dimension of $X_0$. Hence 
 the intersection of $H$ and $X_0$ is a proper closed subset of $X_0$.  Let $T_1$ be the open subscheme $T\backslash H$ of 
 $T$. From the previous discussion, we see that the underlying 
  set of $T_1$ is a {\it non-empty open subset of $X_0$}. 
  
  We are now in a position to complete the proof of Theorem \ref{TVth}. The proof will be by contradiction. 
  So suppose that for all $l\in\mN$, the set 
  $$
\{P\in X_0(\bar k))\ |\ \textrm{$P$ lifts to an element of $\CX(S^\sh_l)\cap p^{l}\cdot\CA(S^\sh_l)$}\}
$$
is Zariski dense in $X_0$. 

Choose an arbitrary $l\in\mN$ and 
 let $P\in T_1(\bar k)$ be a point, which lifts to an element of $\CX(S^\sh_l)\cap p^{l}\cdot\CA(S^\sh_l)$. This exists because 
  the set of points in $X_0(\bar k)$ with this property is assumed to be dense in $X_0$. 
 Let $\widetilde{P}\in\CA(S^\sh_l)$ be such 
  that $p^{l}\cdot\widetilde{P}\in\CX(S^\sh_l)$ and such that $p^{l}\cdot\widetilde{P}_0=P$. Here 
  $\widetilde{P}_0\in A_0(\bar k)$ is the $\bar k$-point induced by $\widetilde{P}$.   Since the map 
  $"p^{l}":\CA_0(\bar k)\to\CA(S^\sh_{l})$ intertwines $\phi$ and 
  $\tau$, we see that 
  $$
  Q(\tau)("p^{l}"(\widetilde{P}_0))=0.
  $$
  By the remark \ref{remQ}, we thus have $$u("p^{l}"(\widetilde{P}_0))\in\TV(S^\sh_{l}).$$ Hence $$"p^{l}"(\widetilde{P}_0)\in T_1(S_{l}^\sh)\subseteq 
  T(S_{l}^\sh).$$ This 
  shows that $T_1(S_{l}^\sh)\not=\emptyset$. Since $l$ was arbitrary, this shows 
  that the generic fibre $T_{1,K}$ of $T_1$ is not empty, which is a contradiction. 
 \endProof

\begin{cor}
We keep the hypotheses of the Theorem \ref{TVth}. 
We suppose furthermore that $\Tor(A(\bar K))\cap X^{+c}_{\bar K}$ is not dense in $X^{+c}_{\bar K}$ for all $c\in \CA(S)$. Then there exists a constant $m\in\mN$, such that for all $c\in\CA(S)$ the set 
$$
\{P\in X^{+c}_0(\bar k))\ |\ \textrm{\rm $P$ lifts to an element of $\CX^{+c}(S^\sh_m)\cap p^{m}\cdot\CA(S^\sh_m)$}\}
$$
is not Zariski dense in $X^{+c}_0$. 
\label{TVuni}
\end{cor}
Here as usual $\CX^{+c}=\CX+c$ is the translate inside $\CA$ of $\CX$ by $c\in\CA(S)$. 
Slightly  abusing notation, we write $X^{+c}$ for $(\CX^{+c})_K$ and 
$X^{+c}_0$ for $(\CX^{+c})_k$. 
\beginProof 
By contradiction. Write $m(\CX^{+c})$ for the smallest integer $m$ such that 
$$
\{P\in X^{+c}_0(\bar k))\ |\ \textrm{$P$ lifts to an element of $\CX^{+c}(S^\sh_{m})\cap p^{m}\cdot
\CA(S^\sh_{m})$}\}
$$
is not Zariski dense in $X_0$.
Suppose that there exists a sequence $(a_n\in\CA(S))_{n\in\mN}$, such 
that $m(\CX^{+a_n})$ strictly increases. Replace $(a_n\in\CA(S))_{n\in\mN}$ by one of its 
subsequences, so that $\lim_{n}a_n=a\in\CA(S)$, where the convergence is for the 
topology given by the discrete valuation on the ring underlying $S$ (notice that 
$\CA(S)$ is compact for this topology, because $S$ is complete and has a finite residue field at its closed point). Replace $(a_n\in\CA(S))_{n\in\mN}$ by one of its subsequences again, so that 
the image of $a_n$ in $\CA(S_n)$ equals the image of $a$ in $\CA(S_n)$. 
By construction, we have $m(\CX^{+a_n})\geqslant n$ and hence 
by definition $m(\CX^{a})\geqslant n$. Since this is true for all $n\geqslant 0$, this contradicts 
Theorem \ref{TVth}.\endProof

The following corollary should be viewed as a curiosity only, since it is a special case of  Theorem \ref{MMprop}. The interest lies in its proof, which avoids the use of jet schemes, unlike the proof of Theorem \ref{MMprop}.

\begin{cor} 
We keep the notations and assumptions of Corollary \ref{TVuni}.   Suppose furthermore that $S$ is a ring of characteristic $p$ and that the fibres of $\CA$ over $S$ are ordinary abelian varieties. 
We also suppose that $\CX$ is smooth over $S$. 
Let $\Gamma\subseteq A(K)$ be a finitely generated subgroup. 
Then the set $X\cap\Gamma$ is not Zariski dense in $X$. 
\label{thord}
\end{cor}

We shall call the topology on $A(K)$ induced by the discrete valuation {\it the $v$-adic topology}. 

Before the proof of the corollary, recall a simple but crucial lemma of Voloch (see \cite[Lemma1]{Voloch-Towards}): 
\begin{lemma}[Voloch]
Let $L_0$ be a field and let $T$ be a reduced scheme of finite type over $L_0$. 
Then $T(L^\sep_0)$ is dense in $T$ if and only if $T$ is geometrically reduced.
\end{lemma}

\beginProof (of Corollary \ref{thord}). The proof is by contradiction. 
We shall exhibit a translate of $X$ by an element of $A(K)$, which violates the conclusion 
of Theorem \ref{TVth}. 
Suppose that $X\cap\Gamma$ is Zariski dense in $X$. 
Let $P_1\in \Gamma$ be such that $(X+P_1)\cap p\cdot \Gamma$ is dense, let $P_2\in p\cdot \Gamma$ such that 
$(X+P_1+P_2)\cap p^2\cdot \Gamma$ is dense in $X$ and so forth. The existence of 
the sequence of point $(P_i)$ is guaranteed by the assumption on $\Gamma$, which implies that 
the group $p^l \Gamma/p^{l+1} \Gamma$ is finite for all $l\geqslant 0$. Since the 
$v$-adic topology on the set $A(K)$ is compact (because $S$ is a discrete valuation ring with a finite residue field),  
the sequence $Q_i=\sum_{l\geqslant 1}^iP_l$ has a subsequence, which converges in $A(K)$. Let $Q$ be the limit point 
of such a subsequence. By construction, $(X+Q)\cap p^l\cdot A(K)$ is dense for all $l\geqslant 0$. 
Let $\CX^{+Q}:=\CX+Q$.

Consider the morphism $([p^l]^*\CX^{+Q})_\red\to\CX^{+Q}$. There is a diagram

\begin{diagram}
([p^l]^*\CX^{+Q})_\red &\rInto & ([p^l]^*\CX^{+Q}) &  &\rTo^{[p^l]}& & \CX^{+Q}\\
                            &          &    \dInto       &  & & &\dInto\\
                            && \CA & \rTo^{\Frob^l_{\CA/S}} & \CA^{(p^l)} & \rTo^{\Ver^l} & \CA
\end{diagram}
where $F_S$ is the absolute Frobenius morphism on $S$, $\CA^{(p^l)}=F^{l,*}_S\CA$ is the base-change of $\CA$ by $F^{l,*}_S$, $\Frob_{\CA/S}$ is the Frobenius morphism relatively to $S$ and 
$\Ver$ is the Verschiebung (see \cite[$\rm VII_A$, 4.3]{SGA3-1} for the latter). The square is cartesian (by definition). 
By assumption, the morphism $\Ver$ is \'etale. Hence $\Ver^{l,*}(\CX^{+Q})$ is a disjoint union of schemes, which are 
integral and smooth over $S$. Let $\CX_1\hookrightarrow\Ver^{l,*}(\CX^{+Q})$ be an irreducible 
component such that $\CX_{1,K}\cap \Frob^l_{A/K}(A(K))$ is dense. Let $\CX_2:=
(\Frob^{l,*}_{\CA/S}(\CX_1))_\red$ be the corresponding reduced irreducible component. 

Now notice that $\CX_{2,K}$ is geometrically reduced, since $\CX_{2}(K)$ is dense 
in $\CX_{2,K}$ (Voloch's lemma). Furthermore $\CX_{2}$ is flat over $S$, because it is reduced and dominates $S$. Hence $(\CX_{2})^{(p^l)}$ is also flat over $S$. Furthermore, by its very construction 
$(\CX_{2,K})^{(p^l)}$ is reduced, since $\CX_{2,K}$ is geometrically reduced. Hence $(\CX_{2})^{(p^l)}$ is reduced (for this last step, see for instance \cite[4.3.8, p. 137]{Liu-Algebraic}). Recall that $(\CX_{2})^{(p^l)}$ stands for the base-change of $\CX_{2}$ by $F^{l,*}_S$. Notice that we have a commutative diagram
\begin{diagram}
\CX_{2} & \rTo^{\Frob^l_{\CX_{2}/S}} & (\CX_{2})^{(p^l)}\\
\dInto &            &   \dInto \\
\CA & \rTo^{\Frob^l_{\CA/S}} & \CA^{(p^l)}
\end{diagram}
and that $\Frob^l_{\CX_{2}/S}$ is bijective. 
Hence $(\CX_{2})^{(p^l)}$ is isomorphic to $\CX_1$. 
Now since $F_S$ is faithfully flat and $\CX_{2}$ is flat over $S$, we see that 
$\CX_{2}$ is actually smooth over $S$, because $\CX_1$ is smooth over $S$.
Hence every point of $\CX_{2}(\bar k)$ can be lifted to a point in 
$\CX_{2}(S^\sh)$ (see for instance \cite[Cor. 6.2.13, p. 224]{Liu-Algebraic}). Since the morphism 
$[p^l]$ is finite and flat and the scheme $\CX^{+Q}$ is integral, we see that the map $\CX_2\to\CX^{+Q}$ is surjective. This implies that the map $\CX_{2}(\bar k)\to\CX^{+Q}(\bar k)$ is surjective. We conclude that 

{\it every element of $\CX^{+Q}(\bar k)$ is liftable to an element in $\CX^{+Q}(S^\sh)\cap p^l\cdot\CA(S^\sh)$.}

Since $l$ was arbitrary, this contradicts Theorem \ref{TVth}. 
\endProof

Now we want to to conclude by 

\begin{remark}\textnormal{In \cite{Buium-Intersections}, A. Buium also introduces an "exceptional set", which is very 
similar to the set $\Exc$ considered here and he makes a similar 
use of it (catching rational points). There is nevertheless one important difference between Buium's and our methods: the proof of 
Theorem \ref{pffth}, which 
is crucial in our study of the structure of $\Exc$ uses  "Galois equations"
and not differential equations as in \cite{Buium-Intersections}. In this sense, our techniques also differ from the techniques employed in \cite{Hrushovski-Mordell-Lang}, which is close in spirit to \cite{Buium-Intersections} and where the galois-theoretic language is not used either.  }
\label{CRa}
\end{remark}

\begin{remark}\textnormal{
Although Corollary \ref{corimp} shows that the Mordell-Lang conjecture may be reduced to 
the Manin-Mumford conjecture under the assumptions of Theorem \ref{MMprop}, the difficulty of circumventing the fact that the underlying 
abelian variety might not be ordinary (which was a hurdle for some some time) is not thus removed. Indeed, the most difficult part of the algebro-geometric proof of the Manin-Mumford conjecture given in \cite{PR2} concerns the analysis of endomorphisms of abelian varieties, 
which are not globally the composition of a separable isogeny with a a power of a relative Frobenius morphism. }
\label{CRb}
\end{remark}

\end{document}